\documentclass[a4paper,12pt]{article}
\usepackage[T2A]{fontenc}
\usepackage[cp1251]{inputenc}
\begin{document}

\author{S.V. Ludkowski.}

\title{Hypercomplex Generalizations of Gaussian-type Measures.}

\date{07 July 2018}
\maketitle

\begin{abstract}
The article is devoted to a new type of measures which are
hypercomplex generalizations of Gaussian-type measures. The
considered such measures are related with solutions of high order
hyperbolic PDEs and related Markov processes. Their characteristic
functionals are investigated. Cylindrical distributions of these
measures are studied. \footnote{key words and phrases: measure;
integral; Markov process; hypercomplex; PDE; hyperbolic 
address: Dep. Appl. Mathematics, Moscow State Techn. Univ. MIREA,
av. Vernadsky 78, Moscow, 119454, Russia \par e-mail:
sludkowski@mail.ru}
\end{abstract}

\section{Introduction.}
Feynman integrals and quasi-measures with values in the field of
complex numbers are very important in mathematics and mathematical
physics, quantum mechanics, quantum field theory and partial
differential equations (PDEs) (see
\cite{dalfom,gulcastb,johnlap,kimfi09,nicolam16} and references
therein). This integral was first introduced by Feynman to present
transition amplitudes for solutions of the Schr\"{o}dinger equation
in quantum mechanics, where the Hamiltonian corresponds to a PDE of
the elliptic type. Since that time they were used for analysis and
solutions of PDEs with elliptic partial differential operators
(PDOs) and also parabolic PDOs with and without the $\bf i$
multiplier. In stochastic analysis measures with values in matrix
algebras or operator algebras on Hilbert spaces are frequently
studied and used for solutions and analysis of PDEs
\cite{dalfom,gihmskorb,gulcastb,meashand,johnlap}.
\par It is necessary to mention that the Feynman integral approach works for
PDOs of order not higher than two, because it is based on measures,
particularly on modifications of Gaussian measures. But if a
characteristic functional $\phi (t)$ of a measure has the form $\phi
(t)=\exp (Q(t))$, where $Q(t)$ is a polynomial, then its degree is
not higher than two as the Marcinkievich theorem states (see, for
example, Ch. II, \S 12 in \cite{shirb11}). \par It is natural to use
hypercomplex numbers for analysis and solutions of PDEs. As it is
known it was Dirac who first used the complexified quaternion
algebra for a solution the Klein-Gordon hyperbolic PDE of the second
order with constant coefficients which is utilized in spin quantum
mechanics \cite{dirac}. On the other hand, for a hyperbolic,
elliptic or parabolic PDO of order higher than two its decomposition
into a product of lower order PDOs is given by Theorems 2.1, 2.7 and
Corollaries 2.2, 2.3, 2.8 in \cite{ludrend14}. Coefficients of
operators in such decompositions are in the corresponding
Cayley-Dickson algebras. Generally the complex field is insufficient
for this purpose.
\par Therefore using this procedure it is possible to reduce a PDE
problem in many cases to a subsequent consideration of PDEs of order
not higher than two with Cayley-Dickson coefficients.
\par Remind that the Cayley-Dickson algebras ${\cal A}_r$ over
the real field $\bf R$ are natural generalizations of the complex
field, where ${\cal A}_2 = \bf H$ denotes the quaternion skew field
(that is the non-commutative field), ${\cal A}_3=\bf O$ denotes the
octonion algebra, ${\cal A}_0=\bf R$, ${\cal A}_1=\bf C$. There are
the canonical embeddings ${\cal A}_r\hookrightarrow {\cal A}_{r+1}$
and each subsequent algebra is obtained by the so called doubling
procedure from the preceding algebra with the help of the doubling
generator \cite{baez,dickson,hamilt,kansol}.
\par Among them the quaternion and octonion algebras
are intensively applied in PDEs, mathematical physics, quantum field
theory, hydrodynamics, industrial and computational mathematics,
non-commutative geometry
\cite{dirac,emch,guetze,guespr,guesprqa,lawmich}. In its turn this
stimulated the development of non-commutative analysis over
quaternions, octonions and the Cayley-Dickson algebras (see
\cite{baez,brdeso,guespr,guesprqa,guetze,luoyst,ludjms7,lufsqv,
lufoclg} and references therein). \par Some results on solutions of
PDEs of order higher than two with the help of Feynman-type
integration over octonions and Cayley-Dickson algebras were shortly
communicated by the author at a conference in \cite{luconrudn17}.
They are based on generalizations of Gaussian measures. Therefore in
order to implement this program of integration of PDEs of order
higher than two it is necessary at first to develop a theory of such
measures.
\par This article is devoted to measures with values in the
complexified Cayley-Dickson algebra. Necessary definitions are
provided. A model PDE leading to such measures is investigated (see
Theorem 2). It appears non only as a generalization of the
corresponding PDE over $\bf R$ or $\bf C$, but also in the process
of decompositions of PDOs of order higher than two into a product of
PDOs of order one or two. Theorems about such measures are proved.
Characteristic functionals of measures and distributions are studied
for this purpose (see Propositions 4, 6 and 7). The first and second
moments of measures are given by Theorem 8. Consistent families of
measures are investigated in Theorem 12. Extensions of norm-bounded
cylindrical Radon distributions are studied in Theorem 14.
\par Within the same subsection $m$ a formula number $n$ is referred as $(n)$
and outside this subsection as m$(n)$.
\par Main results of this work are
obtained for the first time. They open new opportunities for
subsequent studies of Markov processes related with these measures,
PDEs and stochastic PDEs and their solutions including that of
hyperbolic type and parabolic type with hyperbolic and elliptic
terms of orders two or higher, their stochastic processes.

\section{Hypercomplex Measures and PDEs.}

\par {\bf 1.1. Remark. PDEs.} Denote by $\lambda _n$ the Lebesgue measure on the
Euclidean space ${\bf R}^n$. Consider a domain $U$  in ${\bf R}^n$
which is either open, $U=Int (U)$, or canonically closed, $U= cl
(Int (U))$, where $Int (U)$ denotes the interior of $U$, whilst $cl
(U)$ notates the closure of $U$ in ${\bf R}^n$. Recall that the
Sobolev space $H^k(U, \lambda _n,{\cal A}_r)$ is the completion by a
norm $\| f \| _k$ of the space of all $k$ times continuously
differentiable functions $f: U\to {\cal A}_r$ with compact support,
where
$$(1) \quad \| f \| ^2_{k} := \sum_{j=0}^k \int_{U} \| f^{(j)}
(x) \| ^2 \lambda _n(dx),$$ $f^{(j)}(x)=D^j_xf(x)$ denotes the
$j$-th derivative poly-${\bf R}$-linear operator on ${\bf R}^n$ at a
point $x$, where $n$ is a natural number. Particularly, it may be
$U={\bf R}^n$.
\par Suppose that an operator $B_j$ is realized as an elliptic PDO
${\hat B}_j$ of the second order on the Sobolev space $H^{2}({\bf
R}^{m_j}, \lambda _{m_j},{\bf R})$ by real variables $x_{1+\beta _
{j-1}}$,...,$x_{\beta _j}$, where $m_0=0$, $\beta _0=0$, $\beta
_j=m_0+...+m_j$, $~m_j\in \bf N$ for each $j=1, 2,...$.
\par Let $ \{ i_0, i_1,...,i_{2^r-1} \} $  notates the
standard basis of the Cayley-Dickson algebra ${\cal A}_r$ over the
real field $\bf R$ such that $i_0=1$, $i_l^2=-1$ and
$i_li_k=-i_ki_l$ for each $l\ne k$ with $1\le l$ and $1\le k$. Then
${\cal A}_{r,C}$ stands for the complexified Cayley-Dickson algebra
${\cal A}_{r,C} = {\cal A}_{r} \oplus ({\cal A}_{r}{\bf i})$, where
${\bf i}^2=-1$, $~ {\bf i}b=b{\bf i}$ for each $b\in {\cal A}_r$,
$~2\le r<\infty $. Therefore, each complexified Cayley-Dickson
number $z\in {\cal A}_{r,C}$ has the form $z=x+{\bf i}y$ with $x$
and $y$ in ${\cal A}_r$, $x=x_0i_0+x_1i_1+...+x_{2^r-1}i_{2^r-1}$,
while $x_0,...,x_{2^r-1}$ are in $\bf R$. The real part of $z$ is
$Re (z)=x_0=(z+z^*)/2$, the imaginary part of $z$ is defined as $Im
(z)=z-Re (z)$, where the conjugate
 of $z$ is $z^*=\tilde{z} = Re (z) - Im (z)$,
that is $z^*=x^*-{\bf i}y$ with
$x^*=x_0i_0-x_1i_1-...-x_{2^r-1}i_{2^r-1}$. Then
$|z|^2=|x|^2+|y|^2$, where $|x|^2=xx^*=x_0^2+...+x_{2^r-1}^2$. It is
useful also to put $\| z \| = |z| \sqrt{2}$.
\par We consider a second order PDO
of the form $$(2)\quad {\hat B} = - \frac{1}{2} \sum_{j=1}^m a_j
\hat{B}_j,$$ where $a_j=a_{j,0}+{\bf i}a_{j,1}$ are nonzero
coefficients, $~a_j\in {\cal A}_{r,C}$, $~a_{j,0}$ and $a_{j,1}$
belong to ${\cal A}_r$, where ${\hat B}_j$ is an elliptic PDO of the
second order on $H^{2}({\bf R}^{m_j}, \lambda _{m_j},{\bf R})$ by
real variables $x_{1+\beta _{j-1}}$,...,$x_{\beta _j}$,  where $2\le
r <\infty $.
\par There are the natural embeddings $H^{2}({\bf R}^{m_j}, \lambda
_{m_j},{\cal A}_{r,C})\hookrightarrow H^{2}({\bf R}^{n}, \lambda
_{n},{\cal A}_{r,C})$, where $n=m_1+...+m_m=\beta _m$. Thus $\hat B$
and all ${\hat B}_j$ are defined on $H^{2}({\bf R}^{n}, \lambda
_{n},{\cal A}_{r,C})$. Let also $\sigma ^*$ be a first order PDO
$$(3)\quad \sigma ^*f(x)= \sum_{j=1}^m  \sigma _j^*f(x)\mbox{  and}$$
$$(4)\quad \sigma _j^*f(x) = \sum_{k=\beta_{j-1} +1}^{\beta _j}
\psi _{k;j}\frac{\partial  f(x)}{\partial x_k}$$ for each $f \in
H^{1}({\bf R}^{n}, \lambda _{n},{\cal A}_{r,C})$, where $~ \beta _j
=m_0+...+m_j$ for each $j$, $m_0=0$, $\beta _0=0$; $~ \psi _{k;j}\in
{\cal A}_{r,C}$ for each $k$ and $j$. Then the operator
$$(5)\quad {\hat S} = \frac{\partial }{\partial t} + {\hat B}+ \sigma ^* $$
is defined on a Sobolev space $H^{2,1}({\bf R}^n\times {\bf R},
\lambda _{n+1},{\cal A}_{r,C})$, where $H^{k,l}(U\times V, \lambda
_{n+1},{\cal A}_{r,C})$ is the completion relative to a norm $\| f
\|_{k,l}$ of the space of all functions $f(x,t): U\times V\to {\cal
A}_{r,C}$ continuously differentiable $k$ times in $x$ and $l$ times
in $t$ with compact support, where $V$ is either open or canonically
closed in $\bf R$,
$$(6) \quad \| f \| ^2_{k,l} :=
\sum_{j=0}^k \sum_{s=0}^l \int_{U\times V} \| D^j_xD^s_tf(x,t) \| ^2
\lambda _{n+1}(dx),$$ where $x\in U$, $t\in V$. Evidently,
$H^{k,l}(U\times V, \lambda _{n+1},{\cal A}_{r,C})$ has a structure
of a Hilbert space over $\bf R$, also of a two-sided ${\cal
A}_{r,C}$-module. Particularly, $H^{0}(U, \lambda _{n},{\cal
A}_{r,C})= L^2(U, \lambda _{n},{\cal A}_{r,C}).$
\par Using the change of variables we consider operators with constant coefficients
$$(7) \quad {\hat B}_jf(x) = \sum_{u,k=1}^{m_j}
b_{u,k;j} \frac{\partial ^2 f(x)}{\partial x_{u+\beta_{j-1}}
\partial x_{k+\beta_{j-1}}},$$
\\ for each $f \in H^{2}({\bf R}^{n}, \lambda _{n},{\cal A}_{r,C})$,
where $b_{u,k;j}\in {\bf R}$ for every $u, k, j$, $~ \beta _j
=m_0+...+m_j$, $m_0=0$, $\beta _0=0$. We denote by $[B_j]$ a matrix
with matrix elements $b_{u,k;j}\in {\bf R}$ for every $u$ and $k$ in
$\{ 1,...,m_j \} $, where $j=1,...,m$. Also $B_j$ notates a linear
operator $B_j: {\bf R}^{m_j}\to {\bf R}^{m_j}$ prescribed by its
matrix $[B_j]$. Since the operator $\hat{B}_j$ is elliptic, then
without loss of generality the matrix $[B_j]$ is symmetric and
positive definite. Then using a variable change it also is
frequently possible to impose the condition $Re (\psi _{k;j}\psi
_{i;l}^*) =0$ if either $k\ne i$ or $j\ne l$.

\par Let $\sf A$ be a unital normed algebra over
$\bf R$, where $\sf A$ may be nonassociative, let its center $Z({\sf
A})$ contain the real field $\bf R$. Then by $\mbox{}_l\prod_{k=1}^m
u_k$ we denote an ordered product from right to left such that
$$(8)\quad \mbox{}_l\prod_{k=1}^m u_k = u_m(\mbox{}_l\prod_{k=1}^{m-1} u_k)$$
for each $m\ge 2$, where $\mbox{}_l\prod_{k=1}^1 u_k=u_1$;
$u_1,...,u_m$ are elements of $\sf A$. Then we put
$$(9)\quad \exp_l (z) = 1+ \sum_{n=1}^{\infty } \frac{\mbox{ }_l(z^n)}{n!},$$
where $\mbox{ }_l(z^n) =\mbox{}_l\prod_{k=1}^n z$, $~z\in {\sf A}$,
that is for the particular case $u_1=z$,....,$u_n=z$.
\par {\bf 1.2. Definition.} Let $X$ be a right module over ${\cal A}_{r,C}$ such
that \par $X=X_0\oplus X_1i_1\oplus ... \oplus X_{2^r-1}i_{2^r-1}$,
\\ where $X_0$,...,$X_{2^r-1}$ are pairwise isomorphic vector spaces
over $\bf C$. If an addition $x+y$ in $X$ is jointly continuous in
$x$ and $y$ and a right multiplication $xb$ is jointly continuous in
$x\in X$ and $b\in {\cal A}_{r,C}$ and $X_j$ is a topological vector
space for each $j\in \{ 0, 1,...,2^r-1 \} $, then $X$ will be called
a topological right module over ${\cal A}_{r,C}$. \par For the right
module $X$ over ${\cal A}_{r,C}$ an operator $h$ from $X$ into
${\cal A}_{r,C}$ will be called right ${\cal A}_{r,C}$-linear in a
weak sense if and only if it $h(fb)=(h(f))b$ for each $f\in X_0$ and
$b\in {\cal A}_{r,C}$. Then $X^*_r$ denotes a family of all
continuous right ${\cal A}_{r,C}$-linear operators $h: X\to {\cal
A}_{r,C}$ in the weak sense on the topological right module $X$ over
${\cal A}_{r,C}$.
\par An operator $h: X\to {\cal A}_{r,C}$ is called right
${\cal A}_{r,C}$-linear if and only if $h(fb)=(h(f))b$ for each
$f\in X$ and $b\in {\cal A}_{r,C}$.
\par Symmetrically on a left module $Y$ over ${\cal A}_{r,C}$
such that  \par $Y=Y_0\oplus i_1Y_1\oplus ... \oplus
i_{2^r-1}Y_{2^r-1}$,
\\ where $Y_0$,...,$Y_{2^r-1}$ are pairwise isomorphic vector spaces
over $\bf C$ are defined left ${\cal A}_{r,C}$-linear operators and
left ${\cal A}_{r,C}$-linear in a weak sense operators. A family of
all continuous left ${\cal A}_{r,C}$-linear operators $g: Y\to {\cal
A}_{r,C}$ on the topological left module $Y$ over ${\cal A}_{r,C}$
in the weak sense we denote by $Y^*_l$.
\par We say that $X$ is a two-sided module over the complexified
Cayley-Dickson algebra ${\cal A}_{r,C}$ if and only if it is a left
and right module over ${\cal A}_{r,C}$ and $i_jx_j=x_ji_j$ for each
$x_j\in X_j$ and $j\in \{ 0, 1,...,2^r-1 \} $.

\par {\bf 2. Theorem.} {\it Let a PDO $\hat{S}$ be of the form
1$(5)$ fulfilling the  condition \par $(\alpha )$ $Re (a_{j,0})>
|q_j| \cdot |\sin \phi _j|$ with $q_j^2= |Im
(a_{j,0})|^2-|Im(a_{j,1})|^2- 2 {\bf i} Re (a_{j,0}a_{j,1})$,\\
$q_j\in \bf C$, $\phi _j = arg (q_j)$ for each $j$, where $2\le
r<\infty $. Then a fundamental solution ${\cal K}$ of the equation
$$(1)\quad \hat{S}{\cal K}=\delta (x,t) \mbox{  is}$$
$$(2)\quad  {\cal K}(x,t)= \frac{\theta (t)} {(2\pi )^n} \int_{{\bf R}^n} \exp_l ( -
\sum_{j=1}^m \{ \frac{1} {2} a_j (B_j{\bf y}_j,{\bf y}_j) + {\bf
i}({\bf s}_j,{\bf y}_j) \} t ) \exp (-{\bf i}(y,x)) \lambda _n(dy)
,$$ where $\theta (t)=0$ for each $t<0$, while $\theta (t)=1$ for
each $t\ge 0$, where $2\le r<\infty $, $ ~ {\bf y}_j=(y_{\beta
_{j-1}+1},...,y_{\beta _j})$ with $y_k\in {\bf R}$ for each $k$,
$$(3)\quad ({\bf s}_j,{\bf y}_j)=\sum_{k=\beta _{j-1}+1}^{\beta _j} s_ky_k,$$ $
s_{\beta_{j-1}+k}= \psi _{k;j}$ for each $k=1,...,m_j$ and each
$j=1,...,m$.}
\par {\bf Proof.}  To simplify calculations
we consider the Fourier operator $\hat{F}$ and its inverse ${\hat
F}^{-1}$ on $L^2({\bf R}^n,\lambda _n, {\cal A}_{r,C})$:
$$(4)\quad (\hat{F}f)(y) = \int_{{\bf R}^n} f(x) \exp ({\bf i}(y,x))
\lambda _n(dx),$$
$$(5)\quad (\hat{F}^{-1}f)(x) = (2\pi )^{-n} \int_{{\bf R}^n} f(y) \exp (-{\bf i}(y,x))
\lambda _n(dy),$$ where $(y,x) =\sum_{j=1}^ny_jx_j$,
$x=(x_1,...,x_n)$, $x$ and $y$ are in ${\bf R}^n$, $x_k\in {\bf R}$
for each $k=1,....,n$. This implies that
\par $(6)$ $\hat{F}f=\hat{F}_{c}f+(\hat{F}_{s}f){\bf i}$, where
$$(7)\quad (\hat{F}_{c}f)(y) = \int_{{\bf R}^n} f(x) \cos ((y,x))
\lambda _n(dx),$$
$$(8)\quad (\hat{F}_{s}f)(y) = \int_{{\bf R}^n} f(x) \sin ((y,x))
\lambda _n(dx),$$ since by the Euler formula $\exp (\alpha + \beta
{\bf i})=e^{\alpha } (\cos (\beta ) +{\bf i}\sin (\beta ))$ for each
real variables $\alpha $ and $\beta $, also since ${\bf
i}i_j=i_j{\bf i}$ for each $j=0,...,2^r-1$.
\par Then we deduce that $$(4)'\quad \hat{F}(f)=\sum_{j=0}^{2^r-1} \hat{F}(f_j)i_j,
\mbox{  where  }$$ $$(4)''\quad f=\sum_{j=0}^{2^r-1} f_ji_j,
~f_j(x)\in \bf C$$ for each $x\in {\bf R}^n$ and $j=0,...,2^r-1$.
Thus, $\hat{F}^{-1}\hat{F}=\hat{F}\hat{F}^{-1}=I$ on $L^2({\bf
R}^n,\lambda _n, {\cal A}_{r,C})$, where $I$ is the unit operator.
\par We take the $\bf R$-linear space ${\cal S}({\bf R}^n,{\cal A}_{r,C})$
of all infinite differentiable rapidly decreasing at infinity
functions $f: {\bf R}^n\to {\cal A}_{r,C}$. Its topology is provided
by a family of semi-norms:
$$(9)\quad p_{\alpha , \gamma }(f)=\sup_{x\in {\bf R}^n} \| x^{\alpha
}D^{\gamma }_xf(x) \| <\infty $$ for each $\alpha $ in ${\bf N}_0^n$
and $\gamma $ in ${\bf N}_0$, where $x^{\alpha }=x_1^{\alpha
_1}...x_n^{\alpha _n}$, $x=(x_1,...,x_n)$, $x_j\in {\bf R}$ for each
$j$, while ${\bf N} = \{ 1, 2, 3,... \}$, $ ~ {\bf N}_0 = {\bf
N}\cup \{ 0 \} $. Evidently, ${\cal S}({\bf R}^n,{\cal A}_{r,C})=X$
is also the left and right module over the algebra ${\cal A}_{r,C}$.
We consider $X^*_r={\cal S}^*_r({\bf R}^n,{\cal A}_{r,C})$ (see
Definition 1.2). Notice that particularly over the complexified
quaternions if $h\in {\cal S}^*_r({\bf R}^n,{\cal A}_{r,C})$, then
this operator $h$ is right ${\cal A}_{r,C}$-linear operator (see
Definition 1.2), since the algebra ${\bf H}_{\bf C}={\cal A}_{2,C}$
of complexified quaternions is associative. Then we infer that
$(D_{x_j}h;f)= - (h;D_{x_j}f)$ for each $h\in {\cal S}^*_r({\bf
R}^n,{\cal A}_{r,C})$ and $f\in {\cal S}({\bf R}^n,{\cal A}_{r,C})$,
where $(h;f):=h(f)$, since ${\bf i}i_j=i_j{\bf i}$ for each $j\in \{
0, 1,...,2^r-1 \} $.
\par We seek a solution ${\cal K}(x,t)$ which in the $x$ and $t$
variables belongs to ${\cal S}^*_r({\bf R}^n\times {\bf R},{\cal
A}_{r,C})$. Apparently the dual pair $({\cal S};{\cal S}^*_r)$ and
Formulas $(4)'$ and $(4)''$ induce an extension of the Fourier
transform on ${\cal S}^*_r$ by the formula:
$(\hat{F}h;f)=(h;\hat{F}f)$.
\par Applying the Fourier transform
$\hat{F}_{x}$ in the $x$ variable, $x\in {\bf R}^n$, to Equation
$(1)$ and using Formula $(4)$ we infer that
$$(10)\quad \frac{\partial (\hat{F}_{x}{\cal K})(y,t)}{\partial t}
+\frac{1}{2} \sum_{j=1}^m (B_j{\bf y}_j,{\bf y}_j) a_j
((\hat{F}_{x}{\cal K})(y,t)) $$ $$ + \sum_{j=1}^m {\bf i} ({\bf
s}_j,{\bf y}_j) ((\hat{F}_{x}{\cal K})(y,t))= 1(y)\delta (t),$$
where ${\bf y}_j=(y_{1+\beta _{j-1}},..., y_{\beta _j})$; $ ~1(y)=1$
for each $y$, $~ (\delta (t);g(t))=\delta (g)=g(0)$ for each
continuous function $g: {\bf R}\to {\cal A}_{r,C}$, also $(\delta
(x,t); f(x,t))=f(0,0)$ for each continuous mapping $f: {\bf
R}^n\times {\bf R}\to {\cal A}_{r,C}$. \par On the algebra ${\cal
A}_{r,C}$ we take the same norm as in Remark 1.1:
$$(11)\quad \| z\| =\sqrt{2|p|^2+2|q|^2},$$ where $z\in {\cal
A}_{r,C}$, $ ~ z=p+{\bf i}q$, $~p\in {\cal A}_r$ and $q\in {\cal
A}_r$, $|p|=\sqrt{pp^*}$, where $\sqrt{\cdot }$ is the positive
branch of the square root function on $(0, \infty )$. Then
\par $(p+{\bf i}q)(u+{\bf i}v) = (pu-qv)+{\bf i} (pv+qu)$, \\
consequently, $\| zw \|^2= 2|pu-qv|^2+2|pv+qu|^2$ and hence $\| zw
\|^2\le 2(|p||u|+|q||v|)^2+2(|p||v|+|q||u|)^2$ leading to the
estimate
\par $(12)$ $\| zw \| \le \| z \| \| w \|$ \\ for each $z=p+{\bf i}q$ and $w=u+{\bf i}v$ in
${\cal A}_{r,C}$, where $p$, $q$, $u$ and $v$ belong to ${\cal
A}_r$. Next we consider the series
$$(13)\quad \exp_l (z) =1 +\sum_{k=1}^{\infty } \frac{\mbox{ }_l(z^k)}{k!}$$
(see Remark 1.1). It converges on ${\cal A}_{r,C}$, since $$\|
\exp_l(z) \| \le 1+\sum_{k=1}^{\infty } \frac{\| z^k \|}{k!}= \exp
(\| z \| )<\infty .$$ For the exponential function $\exp_l(zt)$ we
deduce that
$$(14)\quad \frac{\partial \exp_l(zt)}{\partial t} = z \exp_l(zt)$$
for each $t\in \bf R$ and $z\in {\cal A}_{r,C}$, since ${\bf
R}\subset \bf C$ and the complex field $\bf C$ is the center
$Z({\cal A}_{r,C})$ of the algebra ${\cal A}_{r,C}$  (see Remark
1.1).

\par Therefore, a solution of the differential equation
$(10)$ is
$$(\hat{F}_{x}{\cal K})(y,t)= \theta (t)\exp_l ( -
\sum_{j=1}^m \{ \frac{1}{2} a_j (B_j{\bf y}_j,{\bf y}_j) + {\bf i}
({\bf s}_j,{\bf y}_j) \} t ).$$ From Formula $(5)$ we deduce Formula
$(2)$, since $(y,x)\in \bf R$ for each $x$ and $y$ in ${\bf R}^n$,
also since $\exp ({\bf i}(y,x))\in {\bf C}$ and ${\bf C}=Z({\cal
A}_{r,C})$ (see Remark 1.1).
\par It remains to verify that this integral converges for each $t>0$, since
$(\hat{F}_{x}{\cal K})(y,t)=0$ for each $t<0$, also $\hat
{F}^{-1}(1)=\delta $. For this mention that
\par $(15)$ $\exp_l (z) = \exp (u_0+{\bf
i}v_0) \exp _l(u'+{\bf i}v')$\\
for each $z\in {\cal A}_{r,C}$, where $u_0$ and $v_0$ are real
numbers, $u'\in {\cal A}_r$, $v'\in {\cal A}_r$, $Re (u')=0$, $Re
(v')=0$ such that $z$ is presented in the form \par $z=u_0+{\bf
i}v_0+u'+{\bf i}v'$. \par We have that $w^2=({u'}^2-{v'}^2)+{\bf i}
(u'v'+v'u')$, where $w=u'+{\bf i}v'$, hence
$w^2=-|u'|^2+|v'|^2+2{\bf i} Re (u'v')$, since $Re (u')=0$, $Re
(v')=0$, \par $(u'v')^*={v'}^*{u'}^*=(-v')(-u')=v'u'$,
$|u'|^2=-{u'}^2$, \\ where $u'$ and $v'$ are in ${\cal A}_r$, $~
2\le r$. This implies that $$\exp_l (w) =\sum_{k=0}^{\infty }
\frac{(w^2)^k}{(2k)!} + w \sum_{k=0}^{\infty }
\frac{(w^2)^k}{(2k+1)!},$$ since $w^2\in \bf C$ and $Z({\cal
A}_{r,C})=\bf C$ (see Remark 1.1). Hence
$$(16)\quad \exp_l (w) = \cos (\nu ) +\frac{w}{\nu } \sin (\nu ),$$
where $\nu \in \bf C$ is such that $w^2=({\bf i}\nu )^2$, extending
Formula $(16)$ also at $\nu =0$ using the fact $\lim_{\nu \to 0} \nu
^{-1}\sin (\nu )=1$. As usually
\par $\cos (\nu ) = \cos (\nu _0) \cdot \cosh (\nu _1) - {\bf i}
\sin (\nu _0) \cdot \sinh (\nu _1)$,
\par $\sin (\nu ) = \sin (\nu _0) \cdot \cosh (\nu _1) + {\bf i}
\cos (\nu _0) \cdot \sinh (\nu _1)$, with \par $\cosh (\nu _1)
=(\exp (\nu _1)+\exp (-\nu _1))/2$ and $\sinh (\nu _1) =(\exp (\nu
_1) - \exp (-\nu _1))/2$, \\ where $\nu =\nu _0+{\bf i}\nu _1$, with
$\nu _0\in \bf R$ and $\nu _1\in \bf R$.
\par Taking $$(17)\quad z(y,t)= - \sum_{j=1}^m \{ \frac{1}{2}
a_j (B_j{\bf y}_j,{\bf y}_j) + {\bf i} ({\bf s}_j, {\bf y}_j) \}
t,$$ we get that \par $u_0 = - \sum_{j=1}^m [\frac{1}{2} (Re
(a_{j,0})) (B_j{\bf y}_j,{\bf y}_j)- Re ({\bf s}_{j,1}, {\bf
y}_j)]t$,
\par $v_0=-\sum_{j=1}^m [Re ({\bf s}_{j,0}, {\bf y}_j) + \frac{1}{2} (Re
(a_{j,1})) (B_j{\bf y}_j,{\bf y}_j)]t$, \par $u' = - \sum_{j=1}^m
[\frac{1}{2} (Im (a_{j,0})) (B_j{\bf y}_j, {\bf y}_j)- Im ({\bf
s}_{j,1}, {\bf y}_j)]t$, \par $v' =- \sum_{j=1}^m [Im ({\bf
s}_{j,0}, {\bf y}_j)+ \frac{1}{2}
(Im (a_{j,1}))(B_j{\bf y}_j,{\bf y}_j)]t$, \\
where ${\bf s}_j={\bf s}_{j,0}+{\bf i} {\bf s}_{j,1}$ with ${\bf
s}_{j,0}$ and ${\bf s}_{j,1}$ belonging to ${\cal A}_r$ for each
$j$, consequently, $u_0=u_0(y,t)$ is a polynomial of the second
order in the variable $y\in {\bf R}^n$, $u'=u'(y,t)$ and
$v_0=v_0(y,t)$ and $v'=v'(y,t)$ are polynomials of order not higher
than two in $y$, since $B_j$ is a positive operator on ${\bf
R}^{m_j}$ and $Re (a_{j,0})>0$ for each $j$. Evidently, $u_0$,
$v_0$, $u'$ and $v'$ are linear in $t$. The operator $B_j$ has the
natural $\bf R$-linear extension $B_j: {\cal A}_{r,C}^{m_j} \to
{\cal A}_{r,C}^{m_j}$. It is ${\cal A}_{r,C}$-linear in the weak
sense: $B_j({\bf y}_jb)= (B_j{\bf y}_j)b=b(B_j{\bf y}_j)$ for each
${\bf y}_j\in {\bf R}^{m_j}$ and $b\in {\cal A}_{r,C}$. The
condition $(\alpha )$ of this theorem and Formulas $(15)$-$(17)$
imply that constants $C_1>0$ and $C_2>0$ exist so that
\par $(18)$ $|u_0(y,t)|^2- |\nu _1(y,t)|^2\ge C_1 |u_0(y,t)|^2$ for each
$|y|>C_2$ and $t\in \bf R$, where $y\in {\bf R}^n$, since \par $\nu
= \sqrt{ |u'|^2-|v'|^2 - 2 {\bf i} Re (u'v')}$, $ ~ \nu _1 =( \nu -
\bar{\nu })/(2{\bf i})$. \par The function $\nu _1(y,t)$ is bounded
on $\{ y\in {\bf R}^n, t\in {\bf R}: |y|\le C_2, |t|\le \rho  \} $
for each $0<\rho <\infty $ and the mapping $z=z(y,t)$ is continuous
on ${\bf R}^n\times \bf R$. Therefore, Formulas $(15)$-$(18)$ imply
that
$$(19)\quad \lim_{|y|\to \infty } |\exp_l (z(y,t))|
\cdot \exp (C_3\sum_{j=1}^mRe (a_{j,0})(B_j{\bf y}_j,{\bf
y}_j)t)=0$$ for each $0<C_3<\frac{1}{2}\sqrt{C_1}$ and $0<t$. On the
other hand, $|\exp (u_0+{\bf i}v_0) | = \exp (u_0)$, consequently,
the integral in Formula $(2)$ absolutely converges for each $t>0$.

\par {\bf 3. Definition.} Let $a_j\in {\cal A}_{r,C}$ satisfy
Condition 2$(\alpha )$ for each $j$, $B_j: {\bf R}^{m_j}\to {\bf
R}^{m_j}$ be a positive definite operator for each $j=1,...,m$,
$p\in {\cal A}_{r,C}^n$, where $n=\beta _m=m_1+...+m_m$, $2\le
r<\infty $. Let also $U: {\cal A}_{r,C}\to {\cal A}_{r,C}$ be an
operator such that $U=\bigoplus_{j=1}^ma_jB_j$, put
$$(1)\quad (y,z)_s= \sum_{k=1}^n y_kz_k$$ for each $y$ and $z$ in
${\cal A}_{r,C}^n$, where $z=(z_1,...,z_n)$, $z_k\in {\cal A}_{r,C}$
for each $k$. Shortly $(y,z)$ will also be written instead of
$(y,z)_s$, when a situation is specified. Then \par $(2)$
$\hat{\vartheta }_{U,p} (y) :=\exp _l(-\frac{1}{2}(Uy,y) + {\bf i}
(p,y) )$ \\ is called a characteristic functional of an ${\cal
A}_{r,C}$-valued measure $\vartheta _{U,p}$ on a Borel $\sigma
$-algebra ${\cal B}({\bf R}^n)$ of the Euclidean space ${\bf R}^n$,
where $y\in {\bf R}^n$. Define a measure $\mu _{U,p}$ on the Borel
$\sigma $-algebra ${\cal B}({\cal A}_{r,C}^n)$ of the two-sided
${\cal A}_{r,C}$-module ${\cal A}_{r,C}^n$ by the formula:
\par $(3)$ $\mu _{U,p}(p+U^{1/2}dh) = \vartheta _{U,p}(dx) \delta
_{U,p} (p'+U^{1/2}dg)$,
\\ where $h=x+g$, $~ h\in {\cal A}_{r,C}^n$, $ ~ x\in {\bf R}^n$,
$~g\in X'$, $~X'=({\cal A}_{r,C}\ominus {\bf R}i_0)^n$, $~{\bf R}^n$
is embedded into ${\cal A}_{r,C}^n$ as ${\bf R}^ni_0$, $p=p_0+p'$
with $p_0\in {\bf R}^n$ and $p'\in X'$,
\par $(4)$ $\int_{Y'} f(y')\delta _{U,p}
(dy')=f(p')$
\\ for each $f\in C^0_b(Y',{\cal A}_{r,C})$, where $C^0_b(Y',{\cal
A}_{r,C})$ denotes the family of all continuous bounded functions
$f$ from $Y'$ into ${\cal A}_{r,C}$, $~Y'=p'+U^{1/2}X'$.

\par {\bf 4. Proposition.} {\it The measure $\mu _{U,p}$ (see Definition 3)
is $\sigma $-additive on ${\cal B}({\cal A}_{r,C}^n)$.}
\par {\bf Proof.} The characteristic functional $\hat{\vartheta }_{U,p}$
(see Formula 3$(2)$) induces operators
$$(1)\quad \vartheta ^R_{U,p}(f) = \int_{{\bf R}^n} \{ \frac{1}{(2\pi )^n}
\int_{{\bf R}^n}\hat{\vartheta }_{U,p}(y)\exp (- {\bf i} (x,y))
\lambda _n(dy) \} f(x)\lambda _n(dx)$$ and
$$(2)\quad \vartheta ^L_{U,p}(f) = \int_{{\bf R}^n}f(x)\{ \frac{1}{(2\pi )^n}
\int_{{\bf R}^n}\hat{\vartheta }_{U,p}(y)\exp (- {\bf i} (x,y))
\lambda _n(dy) \} \lambda _n(dx)$$ for each $f\in L^{\infty } ({\bf
R}^n,\lambda _n , {\cal A}_{r,C})$. Particularly, if $f\in L^{\infty
} ({\bf R}^n,\lambda _n , {\bf R})$, then $\vartheta
^R_{U,p}(f)=\vartheta ^L_{U,p}(f)$, consequently,
\par $(3)$ $\vartheta _{U,p} (Y) = \vartheta ^L_{U,p}(ch_Y)$ \\ is an
${\cal A}_{r,C}$-valued measure on the Borel $\sigma $-algebra
${\cal B}({\bf R}^n)$ of the Euclidean space ${\bf R}^n$, $Y\in
{\cal B}({\bf R}^n)$, where $ch_Y$ denotes the characteristic
function of a set $Y$, where as traditionally $ch _Y(x)=1$ for each
$x\in Y$, also $ch_Y(x)=0$ for each $x\in {\bf R}^n\setminus Y$.
\par From Formulas 2$(15)$, 2$(16)$, 2$(19)$ and 3$(2)$-3$(4)$
it follows that the operator $\mu ^L_{U,p}(f)$ is continuous on
$L^{\infty } ({\cal A}_{r,C}^n,\lambda _{n_r} , {\cal A}_{r,C})$,
where $n_r=2^{r+1}n$, since on $a_j$ condition 2$(\alpha )$ is
imposed for each $j$. Next we consider an integral in which a
function stands on the left to the measure:
$$(4)\quad \mu ^L_{U,p} (f) =\int_Xf(x)\mu _{U,p}(dx)$$ whenever it exists,
where $X={\cal A}_{r,C}^n$.
Evidently,
\par $(5)$ $\mu ^L_{U,p} (bf+cg) = b \mu ^L_{U,p}(f) + c\mu
^L_{U,p}(g)$ \\ for every either ($b$ and $c$ in ${\cal A}_{r,C}$
and $f$ and $g$ in $L^{\infty } ({\cal A}_{r,C}^n,\lambda _{n_r} ,
{\bf C})$) or ($b$ and $c$ in ${\bf C}$ and $f$ and $g$ in
$L^{\infty } ({\cal A}_{r,C}^n,\lambda _{n_r} , {\cal A}_{r,C})$).
\par On the other hand, $ch_Y\in L^{\infty } ({\cal A}_{r,C},
\lambda _{n_r}, {\bf R})$ for each $Y\in {\cal B}({\cal
A}_{r,C}^n)$. Therefore Formula $(5)$ implies that $\mu _{U,p}$ is
the $\sigma $-additive measure on ${\cal B}({\cal A}_{r,C}^n)$,
since ${\bf C}=Z({\cal A}_{r,C})$ (see Remark 1.1).

\par {\bf 5. Corollary.} {\it If conditions of Theorem 2 are fulfilled,
$t>0$, $p=-s$, then the measure $\mu _{Ut,pt}$ is $\sigma $-additive
on ${\cal B}({\cal A}_{r,C}^n)$.}
\par {\bf Proof.} This follows from Proposition 4.

\par {\bf 6. Proposition.} {\it For each $z\in {\cal A}_{r,C}$ the
function $\chi _z(t)=\exp_l(zt)$ is a character from ${\bf R}$
considered as the additive group into the algebra ${\cal A}_{r,C}$
such that
\par $(1)$ $\exp_l (z t) = (\exp_l(z))^t$ and
\par $(2)$ $\chi _z(t+t_1)=\chi _z(t) \chi _z(t_1)$ \\
for each $t\in \bf R$ and $t_1\in \bf R$.}
\par {\bf Proof.} Formula 2$(16)$ applied to $wt$ instead of $w$
gives
$$\exp_l (wt) = \cos (\nu (wt)) +\frac{wt}{\nu (wt)} \sin (\nu (wt)),$$
where $t\in \bf R$, $z=u_0+{\bf i}v_0+w$, $~w=u'+{\bf i}v'$, $Re
(u')=0$, $Re (v')=0$, $u'\in {\cal A}_r$, $v'\in {\cal A}_r$, $u_0$
and $v_0$ are reals, $(wt)^2=({\bf i}\nu )^2$, $\nu $ is a function
of $wt$, $~ \nu =\nu (wt)$, where a branch of the square root
function $\sqrt{ \cdot }$ on $\bf C$ is chosen such that $\sqrt{\xi
}>0$ for each $\xi
>0$. Therefore, $\nu (wt)=|t|\nu (w)$, consequently,
$$(3)\quad \exp_l (wt) = \cos (t \nu (w)) +\frac{w}{\nu (w)} \sin (t\nu (w)),$$
for each $t\in \bf R$, since $\frac{t}{|t|} \sin (|t|\nu (w))=\sin
(t\nu (w))$. This implies that $\exp_l(wt)=(\exp_l(w))^t$ for each
$t\in \bf R$, since $\nu \in \bf C$, $~ w^2\in \bf C$ and
$(wt)^2/\nu ^2(wt)=-1$, $~{\bf C} = Z({\cal A}_{r,C})$. Together
with $\exp ((u_0+{\bf i}v_0)t) = (\exp (u_0+{\bf i}v_0))^t$ and
2$(15)$ which in this situation reads as $\exp_l (zt) = \exp
((u_0+{\bf i}v_0)t) \exp _l(wt)$ this leads to Formula $(1)$ and
hence to $(2)$ also, since $u_0+{\bf i}v_0\in \bf C$.

\par {\bf 7. Proposition.} {\it If $\mu _{U,p}$ is the measure of Definition 3,
then
\par $(1)$ $\mu ^L_{U,p}(id)=p$, where $id(x)=x$ for each
$x\in {\cal A}_{r,C}^n$.}
\par {\bf Proof.} The quadratic form $(B_j{\bf y}_j,{\bf
y}_j)$ on ${\bf R}^{m_j}$ can be extended on ${\cal A}_{r,C}^{m_j}$
by $-(B_j{\bf y}_ji_{2^r}, i_{2^r}{\bf y}_j)$, since
$i_{2^r}c=c^*i_{2^r}$ and $(b i_{2^r})(i_{2^r}c) = - cb$ for each
$b$ and $c$ in ${\cal A}_r$, where $i_{2^r}\in {\cal A}_{r+1}$,
${\cal A}_r\hookrightarrow {\cal A}_{r+1}$ (see also the notation in
\S 2 and \S 3).
\par Form condition 2$(\alpha )$ it follows that the domain
\par $V_{\rho } := \{ w\in {\cal A}_r^{2n}: ~ w=v\oplus y, ~ v\in {\cal A}_r^n,
~ y\in {\cal A}_r^n, ~ z\in {\cal A}_{r,C}^n, ~ z=v+{\bf i}y, ~
\forall ~ j\in \{ 1,...,m \} ~ Re (a_{j,0}(B_j{\bf z}_j, {\bf
z}_j))> |q_j({\bf z}_j)| \cdot |\sin \phi _j({\bf z}_j)|\mbox{ with
}$ \par $(q_j({\bf z}_j))^2= |Im (a_{j,0}(B_j{\bf z}_j, {\bf
z}_j))|^2-|Im(a_{j,1}(B_j{\bf z}_j, {\bf z}_j))|^2$\par $- 2 {\bf i}
Re ([a_{j,0}(B_j{\bf z}_j, {\bf z}_j)]\cdot [a_{j,1}(B_j{\bf z}_j,
{\bf z}_j)])$, $q_j\in \bf C$, $\phi _j({\bf z}_j) = arg (q_j({\bf
z}_j)) \} $\par $ \cup \{ w\in {\cal A}_r^{2n}: ~ |w|< \rho \} $
\\is open in ${\cal A}_r^{2n}$, where $0<\rho <\infty $, $z=({\bf
z}_1,...,{\bf z}_m)$ with ${\bf z}_j\in {\cal A}_{r,C}^{m_j}$ for
each $j$.
\par Using the embeddings ${\bf R}^n\hookrightarrow {\cal A}_r^n
\hookrightarrow {\cal A}_r^{2n}$ we take the ${\cal A}_r$-analytic
extension on the domain $V_{\rho }$ in the $w$-representation of the
characteristic function $\hat{\vartheta }_{U,p}(z)$ with $z=v+{\bf
i}y$ and $w=v\oplus y$ of the measure $\vartheta _{U,p}$. That is
$$(2)\quad \hat{\vartheta }_{U,p}(z)=\int_{{\bf R}^n} \exp ({\bf i}(x,z))
\vartheta ^L_{U,p}(dx)$$ for each $z\in {\bf R}^n$, since $\exp
({\bf i}(x,z)) \in {\bf C}$ for each $x$ and $z$ in ${\bf R}^n$ and
${\bf C}=Z({\cal A}_{r,C})$. The extension of $\hat{\vartheta
}_{U,p}(v+{\bf i}y)$ from ${\bf R}^n$ onto $V_{\rho }$ is specified
using Formulas $(2.2)$, $(2.5)$ and $(2.10)$ in \cite{ludjmsqcf9}
(or see \cite{ludjms7}) and 2$(13)$ and 2$(15)$. In view of Formulas
2$(11)$ and 2$(19)$ the function $\hat{\vartheta }_{U,p}(v\oplus
{\bf i}y)$ is analytic in $v$ and $y$ with $v\oplus y=w\in V_{\rho
}$ and
\par $\lim_{w\in V_{\rho }, ~ |w|\to \infty } \hat{\vartheta }_{U,p}(v+{\bf i}y)=0$ for
each $0<\rho <\infty $. \par Moreover, this implies that
$\hat{\vartheta }_{U,p}(z) = \hat{\vartheta }_{U,p,0}(z) +  {\bf i}
\hat{\vartheta }_{U,p,1}(z)$, where $\hat{\vartheta }_{U,p,0}(z)\in
{\cal A}_r$ and $\hat{\vartheta }_{U,p,1}(z)\in {\cal A}_r$ for each
$v\oplus y=w\in V_{\rho }$ with $z=v+{\bf i}y$, also $\hat{\vartheta
}_{U,p,0}(v+{\bf i}y)$ and $\hat{\vartheta }_{U,p,1}(v+{\bf i}y)$
are ${\cal A}_r$-analytic in $v$ and $y$ with $v\oplus y=w$ in
$V_{\rho }$. Choose $0<\rho <\infty $ such that $({\bf
R}^n-(p_1\oplus p_2))\subset V_{\rho }$, where $p=p_1+{\bf i}p_2$,
$p_1\in {\cal A}_r^n$, $p_2\in {\cal A}_r^n$. Therefore,
$$\mu ^L_{U,p}(id-p)= \int_{{\bf R}^n} (U^{-1/2}(x_0-p_0))
\vartheta _{U,p}(dx_0) $$ $$= \frac{1}{2\pi ^n} \int_{U^{-1/2}{\bf
R}^n} q \{ \int_{U^{-1/2}{\bf R}^n} \exp _l (-{\bf i} (U^{1/2}q+p_0,
U^{1/2}u+p_0)) $$ $$\hat{\vartheta }_{U,p}(U^{1/2}u+p_0)
d(U^{1/2}u+p_0) \} d(U^{1/2}q+p_0),$$ since $\delta
_{U,p}(U^{-1/2}(id(x')-p'))=0$, where $p_0\in {\bf R}^n$, $p'\in
({\cal A}_{r,C}\ominus {\bf R}i_0)^n$, $p=p_0+p'$,
$u=U^{-1/2}(z-p_0)$, $q=U^{-1/2}(x_0-p_0)$, since $\exp _l(\xi )$
coincides with $\exp (\xi )$ on $\bf C$, also since ${\bf
R}^n-p_0={\bf R}^n$. The integral over ${\bf R}^n$ has an extension
to the noncommutative integral with the ${\cal A}_r$-variables
provided by the left algorithm described in
\cite{luoyst,ludjms7,lufsqv}.
\par By virtue of the noncommutative analog of the Jordan Lemma 3.11
in \cite{ludjmsqcf9} (or 2.23 in \cite{ludtsltjms8}) and the change
of variables theorem 2.11 in \cite{ludrfovfejms10} and Formula
2$(19)$ we infer that
$$\int_{{\bf R}^n} (U^{-1/2}(x_0-p_0))
\vartheta _{U,p}(dx_0) = \frac{1}{2\pi ^n} \int_{{\bf R}^n} q \{
\int_{{\bf R}^n} \exp _l (-{\bf i} (q+p_0, u+p_0))\hat{\vartheta
}_{U,p}(u+p_0) du \} dq.$$

\par Using Identities 2$(15)$, 2$(16)$ we deduce that
$$(3)\quad \exp ( - {\bf i} (x,z))
\hat{\vartheta }_{U,p} (z) = \exp_l ( - \frac{1}{2} (Uz,z) - {\bf i}
(x-p,z))$$ for each $x$ and $z$ in ${\bf R}^n$.  \par Hence
$(\hat{F}^{-1}_z\hat{\vartheta }_{U,p}(z))(x) =
(\hat{F}^{-1}_z\hat{\vartheta }_{U,-p}(z))(-x)$ for each $x\in {\bf
R}^n$, where $z\in {\bf R}^n$. Therefore ${\vartheta
}^L_{U,p}(dx)=\vartheta ^L_{U,-p}(-dx)$ for each $x\in {\bf R}^n$.
Making changes of the variables $q\mapsto -q$, $u\mapsto -u$ and
$p\mapsto -p$ in the integral provides $\mu ^L_{U,p}(id-p)=0$.  On
the other hand, $\mu _{U,p}({\cal A}_{r,C}^n)=\hat{\vartheta
}_{U,p}(0)=1$, hence $\mu _{U,p}(id)=p$ due to Formula 4$(5)$.

\par {\bf 8. Theorem.} {\it For the measure $\mu _{U,p}$ of Definition
3
\par $(1)$ $\mu ^L_{U,p}((id_k-p_k)(id_h-p_h)) = a_j b_{k-\beta
_{j-1}, h- \beta _{l-1}; j} \delta _{j,l}$ \\
for each $k$ and $h$ in $ \{ 1,...,n \} $, where $b_{k-\beta _{j-1},
h- \beta _{l-1}; j}\delta _{j,l}$ is a matrix element of the matrix
$[B]=\bigoplus_{j=1}^m [B_j]$ of the operator
$B=\bigoplus_{j=1}^mB_j$ for each $1+\beta _{j-1}\le k \le \beta _j$
and $1+\beta _{l-1}\le h \le \beta _l$, $j=1,...,m$, $l=1,...,m$,
$id_k(x)=x_k$, $x=(x_1,...,x_n)\in {\cal A}_{r,C}^n$, $x_k\in {\cal
A}_{r,C}$.}
\par {\bf Proof.} We take the matrix $[B] = diag ([B_1],...,[B_m])$
written in the block form. We have that each its block $[B_j]$ is
$m_j\times m_j$  a square symmetric positive definite matrix for
each $j$. Other elements outside all blocks $[B_j]$ are zero in
$[B]$. Apparently from the theory of matrices that there are an
orthogonal matrix $[Q]=diag ( [Q_1],...,[Q_m])$ and a diagonal
matrix $[\Lambda ] = diag ([\Lambda _1],...,[\Lambda _m])$ such that
$[Q][B][Q]^t=[\Lambda ]$, where $[Q_j][B_j][Q_j]^t=[\Lambda _j]$ and
$[\Lambda _j] = diag (\lambda _{1+\beta _{j-1}},...,\lambda _{\beta
_j})$ with $\lambda _k>0$ for each $j=1,...,m$ and $k=1,...,n$,
where $[Q][Q]^t=[I]$ and $[Q]^t[Q]=[I]$, where as usually $ ~ [Q]^t$
notates the transposed matrix of $[Q]$, $~ [I]$ is the unit matrix,
since $[B]>0$ (see \cite{gantmb}).
\par Therefore, it is sufficient to consider the case of the
diagonal matrix $[\Lambda ]$ as $[B]$ with
$U=\bigoplus_{j=1}^ma_j\Lambda _j$, where $\Lambda _j$ are operators
corresponding to their matrices $[\Lambda _j]$. From the symmetry
property of the characteristic functional $\hat{\mu }_{U,p}$ it
follows that
\par $\mu ^L_{U,p}(dx'+e_kdx_k)=\mu ^L_{U,{p'}_k-e_kp_k}
(d{x'}_k-e_kdx_k)$ and consequently,
\par $\mu ^L_{U,p}((id_k-p_k)(id_h-p_h)) =0$ for each $k\ne h$, \\
where $e_k=(0,...,0,1,0,...,0)\in {\bf R}^n$ with $1$ at $k$-th
place, ${x'}_k=x-e_kx_k$, ${p'}_k=p-e_kp_k$.
\par Let now $k=l$. Utilizing the identities ${\hat
F}^{-1}_{x'}1=\delta (x')$ and Formulas 3$(3)$, 3$(4)$ we infer,
that
$$\mu ^L_{U,p}((id_k-p_k)^2)=\frac{1}{2\pi } \int_{-\infty }^{\infty
} (x_k-p_k)^2 \{ \int_{-\infty }^{\infty } \exp_l (-a_j\lambda
_ky_k^2 - {\bf i} (x_k-p_k)y_k) dy_k \} dx_k.$$ By virtue of
Proposition 7
$$\mu ^L_{U,p}((id_k-p_k)^2)= \mu ^L_{U,0}(id_k^2)=
\frac{1}{2\pi } \int_{-\infty }^{\infty } z_k^2 \{ \int_{-\infty
}^{\infty } \exp ( - {\bf i} uz) \exp (-\alpha u^2) du \} dz,$$
since $\exp_l(v)=\exp (v)$ for each $v\in {\cal A}_r$, where $\alpha
= a_j\lambda _k$. The function $\exp ( - \alpha u^2)$ is even in the
variable $u$, consequently, $$\mu
^L_{U,p}((id_k-p_k)^2)=\frac{1}{2\pi } \int_{-\infty }^{\infty } z^2
\{ \int_{-\infty }^{\infty } \exp (- \alpha u^2) \cos (uz) du \} dz
= \alpha .$$ This implies Formula $(1)$ with the Kronecker delta
function so that $\delta _{j,l}=0$ for each $j\ne l$, also $\delta
_{j,j}=1$ for each $j$.

\par  {\bf 9. Definitions.} Let $\Omega $ be a set with an algebra
$\cal R$ of its subsets and an ${\cal A}_{r,C}$-valued measure $\mu
: {\cal R}\to {\cal A}_{r,C}$, where $2\le r$, $ ~ \Omega \in \cal
R$. Then
$$(1)\quad |\mu | := \sum_{j=0}^{2^r-1} ( |\mu _{j,0}|+|\mu
_{j,1}| ) $$ is called a variation and $|\mu |(\Omega )$ is a norm
of the measure $\mu $, where $$(2)\quad \mu = \sum_{j=0}^{2^r-1}
(\mu _{j,0}i_j+\mu _{j,1}i_j{\bf i})$$ is the decomposition of the
measure $\mu $, \par $\mu _{j,k}: {\cal R}\to \bf R$, $|\mu _{j,k}|$
denotes the variation of a real-valued measure $\mu _{j,k}$ for each
$j=0,1,...,2^r-1$ and $k=0, 1$,
\par $|\mu |: {\cal R} \to [0, \infty )$.
\par A class $\cal G$ of subsets of a set $\Omega $ is called compact if
for any sequence $G_k$ of its elements fulfilling
$\bigcap_{k=1}^{\infty }G_k=\emptyset $ a natural number $l$ exists
so that $\bigcap_{k=1}^{l}G_k=\emptyset $. \par An ${\cal
A}_{r,C}$-valued measure $\mu $ (not necessarily $\sigma $-additive,
i.e. a premeasure in another terminology) on an algebra $\cal R$ of
subsets of the set $\Omega $ is approximated from below by a class
$\cal H$, where ${\cal H}\subset {\cal R}$, if for each $A\in \cal
R$ and $\epsilon
>0$ a subset $B$ belonging to the class $\cal H$ exists such that
$B\subset A$ and $|\mu |(A\setminus B)<\epsilon $ (see Formula
$(1)$).
\par The ${\cal A}_{r,C}$-valued measure $\mu $ on the algebra $\cal R$
is called Radon, if it is approximated from below by the compact
class $\cal H$. In this case the measure space $(\Omega , {\cal R},
\mu )$ is called Radon.

\par {\bf 10. Remark. Family of ${\cal A}_{r,C}$-valued measures.}
Consider a family of norm-bounded measure spaces $\{ (\Omega ^l,
{\cal R}^l, \mu ^l) : l \in \Lambda \} $, where $\Lambda $ is a set.
Let each $\sigma $-algebra ${\cal R}^l$ be complete relative to a
$\sigma $-additive measure $\mu ^l$. This means by our convention
that ${\cal R}^l={\cal R}^l_{\mu ^l}$, where ${\cal R}^l_{\mu ^l}$
denotes a completion of the $\sigma $-algebra ${\cal R}^l$ relative
to $|\mu ^l|$, where $|\mu ^l|$ is described by Formula 9$(1)$.
\par Put $\Omega _p:=\prod_{l\in \Lambda } \Omega ^l$
to be the product of topological spaces $\Omega ^l$ supplied with
the product, that is Tychonoff, topology $\tau _{\Omega _p}$, where
each $\Omega ^l$ is considered in its $\tau ^l$-topology such that
$\beta (\tau ^l)\subset {\cal R}^l$ and ${\cal R}^l=({\cal B}(\Omega
^l))_{\mu ^l}$ for each $l$, where $\beta (\tau ^l)$ is a base of
the topology $\tau ^l$, $ ~ {\cal B}(\Omega ^l)$ is the Borel
$\sigma $-algebra on $\Omega ^l$ generated by $\tau ^l$. A natural
continuous projection $\pi ^l: \Omega _p\to \Omega ^l$ for each
$l\in \Lambda $ exists. Let ${\cal R}(\Omega _p)$ be the algebra of
the form $\bigcup_{l_1,\dots,l_k\in \Lambda , k\in \bf N} \bigcap
_{m=1}^k\pi _{l_m} ^{-1}({\cal R}^{l_m})$.

\par {\bf 11. Definition. Consistent family of measure spaces.}
Let $\Lambda $ be a directed set (see \cite{eng}) and $\{ (\Omega
^l,{\cal R}^l,\mu ^l ): l\in \Lambda \} $ be a family of
norm-bounded measure spaces, where ${\cal R}^l$ is the algebra (not
necessarily separating). Supply each $\Omega ^l$  with a topology
$\tau ^l$ such that its base $\beta (\tau ^l)$ is contained in a
$\sigma $-algebra ${\cal R}^l$ and ${\cal R}^l=({\cal B}(\Omega
^l))_{\mu ^l}$. Suppose that this family is consistent in the
following sense:
\par $(1)$ a mapping $\pi ^k_l: \Omega ^k\to \Omega ^l$ exists
for each $k\ge l$ in $\Lambda $ such that $(\pi ^k_l)^{-1}({\cal
R}^l) \subset {\cal R}^k$, $\pi ^l_l(x)=x$ for each $x\in \Omega ^l$
and each $l\in \Lambda $, $\pi ^k_l\circ \pi ^m_k=\pi ^m_l$ for each
$m\ge k\ge l$ in $\Lambda $;
\par $(2)$ $\pi ^k_l(\mu ^k)=(\mu ^l)$ for each $k\ge l$ in $\Lambda $.
\par Such family of measure spaces is called consistent.

\par {\bf 12. Theorem.} {\it Suppose that
$\{ (\Omega ^l,{\cal R}^l,\mu ^l ): l\in \Lambda \} $ is a
consistent family of uniformly norm-bounded measure Radon spaces
(see Definition 11), where $\mu ^l: {\cal R}^l\to {\cal A}_{r,C}$
for each $l$, $~2\le r<\infty $ and
\par $(1)$ $\sup_{l\in \Lambda } |\mu ^l|(\Omega ^l)=:{\bf b}
<\infty $. \\ Then a norm-bounded measure space $(\Omega ,{\cal
R}_{\mu },\mu )$ and a mapping $\pi _l: \Omega \to \Omega ^l$ exist
for each $l\in \Lambda $ such that $(\pi _l)^{-1}({\cal R}^l)\subset
{\cal R}_{\mu }$ and $\pi _l(\mu )=\mu ^l$ for each $l\in \Lambda
$.}
\par {\bf Proof.} By the conditions of this theorem the following inclusion
is valid: $(\pi ^k_l)^{-1}({\cal R}^l) \subset {\cal R}^k$ for each
$k\ge l$ in $\Lambda $, hence without loss of generality topologies
$\tau ^k$ can be chosen such that $(\pi ^k_l)^{-1}(\tau ^l) \subset
\tau ^k$ for each $k\ge l$ in $\Lambda $, since each open subset in
$(\Omega ^l,\tau ^l)$ is a union of some subfamily $\cal G$ in
${\cal R}^l$ and $(\pi ^k_l)^{-1}(\bigcup {\cal G})=\bigcup_{A\in
\cal G} (\pi ^k_l)^{-1}(A)$. Therefore, each $\pi ^k_l$ is
continuous and there exists the inverse system ${\sf S}:=\{ \Omega
^k,\pi ^k_l,\Lambda \} $ of the spaces $\Omega ^k$. In view of
Theorem 2.5.2 in \cite{eng} its limit $\lim {\sf S}=:\Omega $ is the
topological space with a topology $\tau _{\Omega }$ and continuous
mappings $\pi _l: \Omega \to \Omega ^l$ such that $\pi ^k_l\circ \pi
_k= \pi _l$ for each $k\ge l$ in $\Lambda $. Moreover, each element
$\omega \in \Omega $ is a thread $\omega =\{ \omega ^l: \omega ^l\in
\Omega ^l$ $\mbox{for each}$ $l\in \Lambda,$ $\pi ^k_l(\omega
^k)=\omega ^l$ $\mbox{for each}$ $k\ge l\in \Lambda \} $. Then $\pi
_l^{-1}({\cal R}^l)=:{\cal G}^l$ is the algebra of subsets in
$\Omega $ for each $l\in \Lambda $. The base of topology of $(\Omega
,\tau _{\Omega })$ consists of subsets $\pi _l^{-1}(A)$, where $A\in
\tau ^l$, $l\in \Lambda $. Since $\Omega ^l$ is supplied with the
topology $\tau ^l$ the base of which is contained in the algebra
${\cal R}^l$ and ${\cal R}^l=({\cal B}(\Omega ^l))_{\mu ^l}$, then
the algebra ${\cal R}:={\cal R}(\Omega ):= \bigcup_{l\in \Lambda
}{\cal G}^l$ generates a $\sigma $-algebra ${\cal R}_{\sigma }$
containing ${\cal B}(\Omega )$.
\par From conditions of 11$(1)$, 11$(2)$ it follows that $
| \mu ^l|(\Omega ^l)\le | \mu ^k|(\Omega ^k)$ for each $l\le k$ in
$\Lambda $, since $\pi ^k_l(\mu ^k)=\mu ^l$ and $(\pi
^k_l)^{-1}({\cal R}^l)\subset {\cal R}^k$. Hence  $\lim_{l\in
\Lambda } | \mu ^l |(\Omega ^l) = \sup_{ l\in \Lambda } | \mu ^l
|(\Omega ^l)$, since the set $\Lambda $ is directed.
\par  On $\cal R$ a function $\mu $ with values
in ${\cal A}_{r,C}$ exists such that $\mu (\pi _l^{-1}(A)):=\mu
^l(A)$ for each $A\in {\cal R}^l$ and each $l\in \Lambda $. If $A$
and $B$ are disjoint elements in $\cal R$, then there exists $l\in
\Lambda $ such that $A$ and $B$ are in ${\cal G}^l$ due to 11$(1)$,
hence
\par  $A=\pi _l^{-1}(C)$ and $B=\pi _l^{-1}(D)$ for some
$C$ and $D$ in ${\cal R}^l$, consequently, $\mu (A\cup B)=\mu
^l(C\cup D)=\mu ^l(C)+\mu ^l(D)=\mu (A)+\mu (B)$, \\ that is, $\mu $
is additive. Then $| \mu | (A)=| \mu ^l| (C)$ for each $A=\pi
_l^{-1}(C)$ with $C\in {\cal R}^l$, hence $| \mu |(\Omega )\le {\bf
b}$ on $\cal R$. This implies that $\mu $ is representable by
Formula 9$(2)$ such that $\mu _{j,k}(\pi _l^{-1}(A))=\mu
_{j,k}^l(A)$ for every $A\in {\cal R}^l$ and $l\in \Lambda $ and
$j=0,...,2^r-1$, $~k=0, 1$.
\par On the other hand, the $\sigma $-additive real-valued
measure $\mu ^l_{j,k}$ has the decomposition $\mu ^l_{j,k} = (\mu
^l_{j,k})^+ - (\mu ^l_{j,k})^-$ with $(\mu ^l_{j,k})^+: {\cal
R}^l\to [0, \infty )$ and $(\mu ^l_{j,k})^-: {\cal R}^l\to [0,
\infty )$ so that $|(\mu ^l_{j,k})^+|\le {\bf b}$ and $|(\mu
^l_{j,k})^-|\le {\bf b}$ for every $l$, $j$ and $k$, since $|(\mu
^l_{j,k})| = (\mu ^l_{j,k})^+ + (\mu ^l_{j,k})^-$.
\par There is the natural embedding of $\Omega $ into the product
$\Omega _p$ (see Proposition 2.5.1 in \cite{eng}).  \par By virtue
of the Prohorov-Kolmogorov theorems (see Theorems IX.4.1 in
\cite{bouimb} and 1.1.4 in \cite{dalfom}) each bounded quasi-measure
$(\mu _{j,k})^{\pm }$ is a $\sigma $-additive measure, consequently,
$\mu _{j,k}$ has an extension to a $\sigma $-additive measure for
each $j$ and $k$. The latter implies by Formula 9$(2)$ that $\mu $
has and extension to a $\sigma $-additive measure on ${\cal R}_{\mu
}$ so that $|\mu |\le {\bf b}$ due to restrictions 11$(1)$, 11$(2)$
and 12$(1)$.

\par {\bf 13. Definition. Cylindrical distribution.} Suppose that $\Omega $ is
a set with an algebra $\cal R$ of it subsets. Suppose also that $\{
(\Omega ,{\cal G}^l,\mu ^l): l\in \Lambda \} $ is a family of
measure spaces such that $\Lambda $ is directed and ${\cal
G}^l\subset {\cal G}^k$ for each $l\le k\in \Lambda $, ${\cal
R}=\bigcup_{ l\in  \Lambda } {\cal G}^l$. Let a mapping $\mu : {\cal
R}\to {\cal A}_{r,C}$ be such that its restriction is \par $(1)$
$\mu |_{{\cal G}^l}=\mu ^l$ and \par $(2)$ $\mu ^k|_{{\cal G}^l}=\mu
^l$ for each $l\le k$ in $\Lambda $. \par Then the triple $(\Omega
,{\cal R},\mu )$ is called the cylindrical distribution. It is
called Radon, if $\mu ^l$ is Radon on ${\cal G}^l$ for each $l\in
\Lambda $.

\par {\bf 14. Theorem.} {\it Suppose that $(\Omega ,{\cal R},\mu )$
is a norm-bounded cylindrical Radon distribution (see Definition
13). Then $\mu $ has an extension to a norm-bounded $\sigma
$-additive measure $\mu $ on a completion ${\cal R}_{\mu }$ of $\cal
R$ relative to $\mu $.}
\par {\bf Proof.} At first we construct a relation in $\Omega $: \\
$x\kappa ^ly$ if and only if for each $S\in {\cal G}^l$ with $x\in
S$ the inclusion $\{ x,y \} \subset S$ is accomplished. Clearly
$x\kappa ^lx$, that is, $\kappa ^l$ is reflexive. Apparently the
relation $x\kappa ^ly$ is equivalent to $y\in V_x:=\bigcap_{x\in
S\in {\cal G}^l} S$. Then from $y\in S\in {\cal G}^l$ and $y\in V_x$
it follows, that $x\in S$, since otherwise $y\notin V_x$, because
${\cal G}^l$ is a $\sigma $-algebra, $\Omega \in {\cal G}^l$.
Therefore, $V_x=V_y$ and $y\kappa ^lx$, hence the relation $\kappa
^l$ is symmetric. If $x\kappa ^ly$ and $y\kappa ^lz$, then
$V_x=V_y=V_z$, hence $x\kappa ^lz$. The latter means that the
relation $\kappa ^l$ is transitive. Therefore, $\kappa ^l$ is the
equivalence relation (see also Section 2.4 in \cite{eng}). By virtue
of Proposition 2.4.3 in \cite{eng} there exists a quotient mapping
$\pi _l: \Omega \to \Omega ^l:=\Omega /\kappa ^l$, where $\Omega ^l$
is supplied with the quotient topology $\tau ^l$. This generates the
$\sigma $-algebra ${\cal R}^l={\cal G}^l/\kappa ^l$ such that $\beta
(\tau ^l)\subset {\cal R}^l$ and $({\cal B}(\Omega ^l))_{\bar{\mu
}^l} = {\cal R}^l$, where $\bar{\mu }^l=\pi _l(\mu ^l)$ is the
projection of the measure $\mu ^l$ from ${\cal G}^l$ onto ${\cal
R}^l$. This implies that $\pi ^{-1}_l({\cal R}^l)= {\cal G}^l$ and
$\pi _l({\cal G}^l)={\cal R}^l$. \par On the other hand, ${\cal
G}^k\supset {\cal G}^l$ for each $k\ge l$ in $\Lambda $, hence the
topology $\tau ^l$ can be chosen for each $l\in \Lambda $ such that
on $(\Omega ^k,(\pi ^k_l)^{-1}({\cal R}^l))$ an equivalence relation
$\kappa ^k_l$ and a quotient (continuous) mapping $\pi ^k_l: \Omega
^k\to \Omega ^l$ exists such that $\pi ^k_l\circ \pi ^m_k=\pi ^m_l$
for each $l\le k\le m$ in $\Lambda $. Thus an inverse mapping system
$\{ \Omega _k,\pi ^k_l,\Lambda \} $ exists. Therefore, the set
$\Omega $ in the weakest topology $\tau _{\Omega }$ relative to
which the mapping $\pi _l: \Omega \to \Omega ^l$ is continuous for
each $l\in \Lambda $ is homeomorphic with $\lim \{ \Omega ^k,\pi
^k_l,\Lambda \} $ by Proposition 2.5.5 in \cite{eng}.
\par From the above construction it follows that
$\pi _l(\mu )|_{{\cal R}^l}=\bar{\mu }^l$ is a bounded measure on
$(\Omega ^l,{\cal R}^l)$ such that $\pi ^k_l(\bar{\mu }^k)=\bar{\mu
}^l$ and $(\pi ^k_l)^{-1}({\cal R}^l)\subset {\cal R}^k$ for each
$k\ge l\in \Lambda $. Therefore, $\{ (\Omega ^l,{\cal R}^l,\mu ^l):
l\in \Lambda \} $ is the consistent uniformly norm-bounded family of
measure spaces, since $\mu $ is bounded. Thus the statement of this
theorem follows from that of Theorem 12.

\par {\bf 15. Remark. Product of measure spaces.}
Suppose that $\Omega :=\prod_{t\in T}\Omega _t$ is a product of sets
$\Omega _t$ and on $\Omega $ an algebra $\cal R$ of its subsets is
given so that for each $n\in \bf N$ and pairwise distinct points
$t_1<\dots <t_n$ in a set $T\subset {\bf R}$ an ${\cal
A}_{r,C}$-valued measure $P_{t_1,\dots,t_n}$ exists on an algebra
${\cal R}_{t_1,\dots,t_n}$ of $\Omega _{t_n}\times \cdots \times
\Omega _{t_1}$ such that $\pi ^{t_1,\dots,t_n}_{t_1,\dots
,t_{j-1},t_{j+1},\dots ,t_n}({\cal R}_{t_1,\dots,t_n})={\cal
R}_{t_1,\dots ,t_{j-1},t_{j+1},\dots , t_n}$ for each $1\le j\le n$
and $P_{t_1,\dots,t_n}(A_n\times \cdots \times A_{j+1}\times \Omega
_j\times A_{j-1}\times \cdots \times A_1)= P_{t_1,\dots
,t_{j-1},t_{j+1},\dots ,t_n}(A_n\times \cdots \times A_{j+1}\times
A_{j-1}\times \cdots \times A_1)$ for each $A_n\times \cdots \times
A_{j+1}\times  A_{j-1}\times \cdots \times A_1\in {\cal
R}_{t_1,\dots ,t_{j-1},t_{j+1},\dots ,t_n}$, where $\pi
^{t_1,\dots,t_n}_{t_1,\dots ,t_{j-1},t_{j+1},\dots ,t_n}: \Omega
_{t_n}\times \cdots \times \Omega _{t_1}\to  \Omega _{t_n}\times
\cdots \times \Omega _{j+1}\times \Omega _{j-1}\times \cdots \times
\Omega _{t_1}$ is the natural projection, $A_l\subset \Omega _{t_l}$
for each $l=1,\dots,n$. Let this cylindrical distribution be
norm-bounded:
\par $\sup_{t_1<\dots <t_n\in T, n\in \bf N}
| P_{t_1,\dots,t_n}|(\Omega _{t_n}\times \cdots \times \Omega
_{t_1}) <\infty $ and the limit exists
\par $\lim_{t_1<\dots <t_n \in T_0; n\in \bf N}
P_{t_1,\dots,t_n}(\Omega _{t_n}\times \cdots \times \Omega
_{t_1})\in {\cal A}_{r,C}$, where $T_0:=\{ t\in T: P_t(\Omega _t)\ne
0 \} $, $T\setminus T_0$ is finite.
\par Apparently the product of measure spaces is a particular case of
a cylindrical distribution.
\par Therefore, the following corollary evidently follows from this remark and
Theorem 14.
\par {\bf 16. Corollary.} {\it The cylindrical Radon distribution
$P_{t_1,\dots,t_n}$ satisfying conditions of Remark 15 has an
extension to a norm-bounded measure $P$ on a completion ${\cal R}_P$
of ${\cal R}:=\bigcup_{t_1<\dots <t_n\in T, n\in \bf N} {\cal
G}_{t_1,\dots,t_n}$ relative to $P$, where ${\cal
G}_{t_1,\dots,t_n}:= (\pi _{t_1,\dots,t_n})^{-1}({\cal
R}_{t_1,\dots,t_n})$ and $\pi _{t_1,\dots,t_n}: \Omega \to \Omega
_{t_n}\times \cdots \times \Omega _{t_1}$ is the natural
projection.}


\begin{thebibliography}{99}
\bibitem{baez} J.C. Baez.  Bull. Amer.
Mathem. Soc. "The octonions", {\bf 39: 2} (2002), 145-205.

\bibitem{bouimb} N. Bourbaki. "Int\'egration". Livre VI.
Fasc. XIII, XXI, XXIX, XXXV. Ch. 1--9. (Paris: Hermann, 1965, 1967,
1963, 1969).

\bibitem{brdeso} F. Brackx, R. Delanghe, F. Sommen.
"Clifford analysis" (London: Pitman, 1982).

\bibitem{dalfom} Yu.L. Dalecky, S.V. Fomin. "Measures and
differential equations in infinite-dimensional space" (Dordrecht:
Kluwer, 1991).

\bibitem{dickson} L.E. Dickson. "The collected mathematical papers".
Volumes 1-5 (New York: Chelsea Publishing Co., 1975).

\bibitem{dirac} P.A.M. Dirac. "Die Prinzipen der Quantenmechanik"
(Leipzig: Hirzel, 1930).

\bibitem{emch} G. Emch. Helv. Phys. Acta.
"M$\grave e$chanique quantique quaternionienne et
Relativit$\grave e$ restreinte", {\bf 36} (1963), 739-788.

\bibitem{eng} R. Engelking. "General topology" (Moscow: Mir, 1986).

\bibitem{gantmb} F.R. Gantmacher. "Theory of matrices"
(Moscow: Nauka, 1988).

\bibitem{gihmskorb} I.I. Gihman, A.V. Skorohod. "The
theory of stochastic processes" (New York: Springer-Verlag, 1975).

\bibitem{gulcastb} A. Gulisashvili, J.A. van Casteren.
"Non-autonomous Kato classes and Feynman-Kac propagators" (New
Jersey: World Scientific, 2006).

\bibitem{guespr} K. G\"urlebeck, W. Spr\"ossig. "Quaternionic and
Clifford calculus for physicists and engineers" (Chichester: John
Wiley and Sons, 1997).

\bibitem{guesprqa} K. G\"urlebeck, W. Spr\"ossig. "Quaternionic
analysis and elliptic boundary value problem" (Basel: Birkh\"auser,
1990).

\bibitem{guetze} F. G\"ursey, C.-H. Tze. "On the role of
division, Jordan and related algebras in particle physics"
(Singapore: World Scientific Publ. Co., 1996).


\bibitem{hamilt} W.R. Hamilton. "Selected papers. Optics. Dynamics.
Quaternions" (Moscow: Nauka, 1994).

\bibitem{meashand} "Handbook of measure theory", Edited by E. Pap.
V. 1, 2. (Amsterdam: Elsevier, 2002).

\bibitem{johnlap} G.W. Johnson, M.L. Lapidus. "The Feynman integral
and Feynman's operational calculus" (New York: Clarendon Press,
Oxford Univ. Press, 2000).

\bibitem{kansol} I.L. Kantor, A.S. Solodovnikov.
"Hypercomplex numbers" (Berlin: Springer-Verlag, 1989).

\bibitem{kimfi09} B.J. Kim, B.S. Kim. Integral
Transforms and Special Functions. "Integration by parts formulas for
analytic Feynman integrals of unbounded functionals", {\bf 20: 1}
(2009), 45-57.

\bibitem{lawmich} H.B. Lawson, M.-L. Michelson. "Spin geometry"
(Princeton: Princ. Univ. Press, 1989).

\bibitem{ludmmas2010} S.V. Ludkovsky. Mathemat. Methods in the Appl. Sci. "Feynman integration
over octonions with application to quantum mechanics", {\bf 33: 9}
(2010); 1148-1173.

\bibitem{ludjmsqcf9} S.V. Ludkovsky. J. Mathem. Sci., N.Y. "Noncommutative quasi-conformal
integral transforms over quaternions and octonions", {\bf 157: 2}
(2009), 199-251.


\bibitem{luoyst} S.V. Ludkovsky, F. van Oystaeyen. Bull. Sci. Math.
(Paris). Ser. 2. "Differentiable functions of quaternion variables",
{\bf 127} (2003), 755-796.

\bibitem{ludjms7} S.V. Ludkovsky.  J. Mathem. Sci., N.Y.
"Differentiable functions of Cayley-Dickson numbers and line
integration", {\bf 141: 3} (2007), 1231-1298.

\bibitem{lufsqv} S.V. Ludkovsky. J. Mathem. Sci., N.Y.
"Functions of several Cayley-Dickson variables and manifolds over
them", {\bf 141: 3} (2007), 1299-1330.

\bibitem{lufoclg} S.V. Ludkovsky.
"Geometric loop groups and diffeomorphism groups of manifolds,
stochastic processes on them, associated unitary representations".
In the book: "Focus on Groups Theory Research" (Nova Science
Publishers, Inc.: New York) 2006, pages 59-136.

\bibitem{ludtsltjms8} S.V. Ludkovsky. J. Mathem. Sci., N.Y.
"Two-sided Laplace transform over Cayley-Dickson algebras and its
applications", {\bf 151: 5} (2008), 3372-3430.

\bibitem{ludrfovfejms10} S.V. Ludkovsky. Far East J. of Math.
Sci. (FJMS). "Residues of functions of octonion varaibles", {\bf 39:
1} (2010), 65-104.

\bibitem{ludrend14} S.V. Ludkowski. Rendic. dell'Ist. di Math. dell'Universit\`a di Trieste. Nuova
Serie; "Decompositions of PDE over Cayley-Dickson algebras", {\bf
46} (2014), 1-23.

\bibitem{ludcmft12} S.V. Ludkovsky. Comput. Methods and Function Theory.
"Line integration of Dirac operators over octonions and
Cayley-Dickson algebras", {\bf 12: 1} (2012), 279-306.

\bibitem{luconrudn17} S.V. Ludkovsky. Feynman-type local integration
of stochastic PDE. P. 453-455. In: Analytical and computational
methods in probability theory and its applications (ACMPT-2017),
Proceedings of the International Scientific Conference.  23-27
October 2017, Moscow, Russia, Ed. A.V. Lebedev.

\bibitem{nicolam16} F. Nicola. Advances in Mathematics. "Convergence in $L^p$ for Feynman
path integrals", {\bf 294} (2016), 384-409.

\bibitem{shirb11} A.N. Shiryaev. "Probability" (Moscow: MTsNMO,
2011).


\end{thebibliography}
\end{document}